\documentclass[preprint,12pt]{elsarticle}
\usepackage{fancyhdr}
\usepackage[T1]{fontenc}
\usepackage[latin1]{inputenc}
\usepackage[english]{babel}
\usepackage{times}
\usepackage{a4, amsfonts, latexsym, amscd, 
  amsmath, amsthm,
  amssymb,
  enumerate,
  bbm,
  textcomp,
  appendix,
  tabularx}
  \usepackage{float}
\usepackage{tikz}
\usepackage{tikz-cd}
\usepackage[all]{xy}
\usepackage[active]{srcltx}
\usepackage{hyperref}
\usepackage{verbatim}
\usepackage{color}
\usepackage{setspace}
\usetikzlibrary{decorations.pathreplacing}

\newtheorem*{theo}{Theorem}
\newtheorem*{lem}{Lemma}
\newtheorem*{prop}{Proposition}
\newtheorem*{cor}{Corollary}
\newtheorem*{defn}{Definition}

\newtheorem*{ex}{Example}

\newcommand{\add}{\mathrm{add}}

\newcommand{\gen}{\mathrm{Gen}}
\newcommand{\cogen}{\mathrm{Cogen}}
\newcommand{\Mod}{\mathrm{mod}}
\newcommand{\md}{\mathrm{mod}}
\newcommand{\ind}{\mathrm{ind}}
\newcommand{\Hom}{\mathrm{Hom}}

\newcommand{\Ext}{\mathrm{Ext}}
\newcommand{\End}{\mathrm{End}}

\newcommand{\pd}{\mathrm{pd}}
\newcommand{\id}{\mathrm{id}}

\newcommand{\La}{\mathcal{L}}

\newcommand{\theend}{\hspace*{\fill} $\square$\\}
\newcommand{\Ker}{\mathrm{Ker}}

\frontmatter

 \title{On the global dimension three endomorphism algebras of the minimal generator-cogenerator}

\author{Edson Ribeiro Alvares}

\address{Centro Polit\'ecnico, Departamento de Matem\'atica,
  Universidade Federal do Paran\'a, CP019081, Jardim das Am\'ericas,
  Curitiba-PR, 81531-990, Brazil }

\ead{rolo1rolo@gmail.com, rolo@ufpr.br}

\author{Clezio Aparecido Braga}

\address{Centro de Ci\^encias Exatas e Tecnol\'ogicas, Universidade
  Estadual do Oeste do Paran\'a, Cascavel-PR, 85819-110, Brazil}
\ead{clezio.braga@unioeste.br}

\author{Sonia Trepode}

\address{CEMIM, FCEyN, Universidad Nacional de Mar del Plata.CONICET. Argentina}

\ead{strepode@gmail.com}

\author{Heily Wagner}

\address{Centro Polit\'ecnico, Departamento de Matem\'atica,
  Universidade Federal do Paran\'a, CP019081, Jardim das Americas,
  Curitiba-PR, 81531-990, Brazil }

\ead{heilywagner@ufpr.br}

\begin{document}

 \begin{abstract}
   The main goal of this paper is to study the class of algebras for which the global dimension of the endomorphism ring of the generator-cogenerator, given by the sum of the projective and injective modules, is equal to three. We will refer to these algebras as representation-hereditary algebras. We show that these algebras are torsionless-finite, as defined by Ringel. These algebras do not necessarily have finite global dimension; however, when there is no non-zero morphism from an injective to a projective module, they have global dimension less than or equal to two, with some additional homological properties. By utilizing the general framework provided by the study of the representation dimension of an algebra, we present further homological consequences. In the case where these algebras are tame quasi-tilted algebras, we prove that they belong to  certain classes of tilted algebras. Although not all tilted algebras are representation-hereditary, we provide sufficient conditions for them to be representation-hereditary.
\end{abstract}

\maketitle

\section{Introduction}

The endomorphism ring of generator-cogenerator modules has been
extensively studied (see \cite{MR349740}, \cite{MR2494382},
\cite{MR3078656}). A specific study on the case of hereditary algebras
was conducted by Dlab and Ringel in \cite{MR2494382}.  A particularly
influential perspective on this subject was developed by Auslander,
who clearly emphasized the importance of the global dimension of such
endomorphism rings through what became known as the \emph{Auslander
correspondence} (see \cite{MR1674397}, \cite{MR1476671})).

In this context, Auslander defined the \emph{representation dimension}
of an Artin algebra $A$ as the infimum of the global dimensions of
endomorphism rings of generator-cogenerators in $\mathrm{mod}\, A$;
that is, $\mathrm{rep.dim}(A) := \inf \{ \mathrm{gl.dim}\,
\mathrm{End}_A(M) \mid M \text{ is a generator-cogenerator in }
\mathrm{mod}\, A \}$.

It is known that if $A$ is semisimple, then its representation
dimension is zero. Moreover, $A$ has representation dimension equal to
two if and only if it is of finite representation type. A
generator-cogenerator module $M$ for which $\mathrm{gl.dim}\,
\mathrm{End}_A(M)$ attains the representation dimension is called an
\emph{Auslander generator}. In the case of hereditary algebras of
infinite representation type, it is known that the Auslander generator
is $A \oplus DA$ and that $\mathrm{gl.dim}\, \mathrm{End}_A(A \oplus
DA) = 3$.

The aim of the present work is to study the class of algebras for which $\mathrm{gl.dim}, \mathrm{End}_A(A \oplus DA) = 3$. We refer to these algebras as \emph{representation-hereditary algebras}. Note that not every hereditary algebra is representation-hereditary. For instance, if we consider the algebras given by the oriented quivers $Q = \vec{A}_2$ or $\vec{A}_3$, we obtain $\mathrm{gl.dim} \mathrm{End}_A(A \oplus DA) = 2$ and $3$, respectively. However, when $A$ is a hereditary algebra of infinite representation type, this dimension is always equal to three.

We have two reasons for naming these algebras representation-hereditary. One reason is that they contain all hereditary algebras of infinite representation type. However, the main reason arises from the homological aspects presented in the next theorem, namely, the kernel of the right approximation is projective, and the cokernel of the left approximation is injective. The results we obtain highlights some homological aspects of
representation-hereditary algebras.

This result stems from the
following fact: a method for computing the representation dimension is
through approximation resolutions. 
Therefore, the following result describes precisely what these
approximation resolutions are in the specific case where the
generator-cogenerator is $A \oplus DA$.

\begin{theo}[Theorem \eqref{kernel_projective}]
For an artin algebra $A$ such that $\Mod A \ne \add (A \oplus DA)$ the
following are equivalent:
\begin{enumerate}
    \item $A$ is a representation-hereditary algebra.
    \item for all $M \in \ind A \setminus \add (A \oplus DA)$ there
      exists a right $\add (A \oplus DA) $-approximation $f \colon P_0
      \oplus I_0 \to M$ such that $Ker f$ is projective.
    \item for all $M \in \ind A \setminus \add (A \oplus DA)$ the
      kernel of a minimal right $\add (A \oplus DA) $-approximation $f
      \colon P_0 \oplus I_0 \to M$ is projective.
    \item for all $M \in \ind A \setminus \add (A \oplus DA)$ there
      exists a left $\add (A \oplus DA) $-approximation $f \colon M
      \to P^0 \oplus I^0$ such that $Coker f$ is injective.
    \item for all $M \in \ind A \setminus \add (A \oplus DA)$ the
      cokernel of a minimal left $\add (A \oplus DA) $-approximation
      $f \colon M \to P^0 \oplus I^0$ is injective.
\end{enumerate}
\end{theo}

Another relevant aspect of this class of algebras is the fact that the
indecomposable modules in $\mathrm{Gen}\, DA$ and those in
$\mathrm{Cogen}\, A$ are, respectively, cosyzygies of indecomposable
projective modules and syzygies of indecomposable injective modules.
This property places this class of algebras within the broader class
known as \emph{torsionless-finite algebras} (see
\cite{zbMATH06154188}).

\begin{theo}[Theorem \eqref{torsionless}]
  Let $A$ be a representation-hereditary algebra. Then
  \begin{enumerate}
  \item[(a)] If $M$ is a non injective indecomposable module in $\gen
    DA$ then $M = \Omega^{-}(P) $ for some projective indecomposable
    module $P$.
     \item[(b)] If $M$ is a non projective indecomposable module in
       $\cogen A$ then $M = \Omega(I) $ for some injective
       indecomposable module $I$.
\item[(c)] $A$ is torsionless-finite.
\end{enumerate} 
\end{theo}

While every representation-hereditary algebra is torsionless-finite,
not every torsionless-finite algebra is representation-hereditary. It
suffices to consider the hereditary algebra we discussed earlier. This
issue highlighted the need to search for sufficient conditions under
which an algebra is representation-hereditary. The intricate behavior
of the module category of this class of algebras made it possible to
identify such sufficient conditions.  These conditions follow directly from the
homological assumptions of Theorem~\eqref{kernel_projective}.

The sufficient conditions are as follows:
\begin{theo} [Theorem \eqref{suficiente}]
Let $A$ be an Artin algebra satisfying the following properties:
\begin{enumerate}
\item[(b.1$ $)] $ {\rm ind}A \cap \Ker(\mathrm{Hom}(\_, A)) \subseteq
  \mathcal{I}^{\leq 1}$;
\item[(b.2)] $\mathrm{ind}A \setminus \mathcal{P}^{\leq 1} \subseteq
  \mathrm{add}DA;$
\item[(b.3)] $\mathrm{Supp}(\mathrm{Hom}(\_, A)) \setminus
  \mathrm{add}DA \subseteq \mathrm{Cogen}A$;
\end{enumerate}
Then $A$ is a representation-hereditary algebra.
\end{theo}

Another relevant aspect of this class of algebras is the fact that there is no upper bound for the global dimension. However, if $\mathrm{Hom}_A(DA, A) = 0$, then the global dimension is at most two. Within this context - namely, the context of algebras with such an upper bound on global dimension - we have the class of quasi-tilted algebras. Not every quasi-tilted algebra is representation-hereditary. In general, quasi-tilted algebras that are not tilted satisfy $\mathrm{gl.dim}, \mathrm{End}_A(A \oplus DA) \leq 4$ (see \cite[Proposition 3.1]{APT06}).

Therefore, under this condition, representation-hereditary algebras that are not quasi-tilted exhibit a distinctive behavior concerning the orbits of their modules with projective and injective dimension equal to two.

\begin{theo}[Theorem \eqref{quasitilted}]
  Let $A$ be a representation-hereditary algebra. If $\Hom(DA, A) = 0$, then one of the following conditions hold:
  \begin{enumerate}
  \item[(a)] $A$ is a quasitilted algebra or
   \item[(b)] There exist an indecomposible module $X$ such that  $\mathrm{pd}\, X = \mathrm{id}\, X = 2$, then $  \mathrm{pd}\, \tau^{-1} X = \mathrm{id}\, \tau X = 2.$
\end{enumerate}
  \end{theo}

Not every quasitilted algebra is representation-hereditary. Nevertheless, the following result allows us to describe, in the tame case, the classes of algebras in which the representation-hereditary algebras are contained. In the proof of this result, we make essential use of the fact that a certain type of component of the Auslander-Reiten quiver does not occur as a component in representation-hereditary algebras.

\begin{theo}[Theorem \eqref{concealed}]
  Let $A$ be a representation-hereditary algebra.  
  If $A$ is a tame quasi-tilted algebra, then $A$ is tilted of one of the following types:
  \begin{enumerate}
   \item[(a)] $A$ is tilted of finite representation type.
   \item[(b)] $A$ is tame concealed.
   \item[(c)] $A$ is tame tilted of wild type, without inserted or coinserted tubes, and the slice of $A$ lies in the connecting component, which is neither a postprojective nor a preinjective component of $A$.
  \end{enumerate}
\end{theo}

In the particular case of tilted algebras, the necessary conditions (Theorema \eqref{suficiente}) admit a clear interpretation. For instance, hypothesis $(b_{2})$ can be reformulated in terms of $\mathrm{Cogen}(\tau_H T)$ (see  \eqref{cogen_tau}). These interpretations enabled the derivation of sufficient conditions for a tilted algebra to be representation-hereditary:

\begin{theo}[Theorem \eqref{theo_tilted}]
Let $A$ be a tilted algebra, $H$ be a hereditary algebra and $T_H$ be a tilting module such that $A = \End T$. Denote by $S\subseteq (Q_H)_0$ the set of all sink in $Q_H$, $R = \{i \in (Q_H)_0 \mid P_H(i) \in \add T\}$  and let $I = \bigoplus_{r \in R} I_H(r)$. If
\begin{enumerate}
    \item[(1)] $S \subseteq R$ ;
    \item[(2)] $\cogen \tau_H T = \add \tau_H T$;
    \item[(3)]  for any  $i \notin R$, if $\varphi \colon T_o \oplus I_0 \to I_H(i)$ is a right minimal $\add(T \oplus I)$-approximation of $I_H(i)$ then $\mathrm{Ker} \varphi \in \add T$.
\end{enumerate}
then $A$ is a representation-hereditary algebra.
\end{theo}

This paper is organized as follows: 

In Section~2, we present a summary of the notation and briefly recall
the basic elements of the theory. In Section~3, we introduce the class
of representation-hereditary algebras and present one of our main
results from a homological perspective.  Section~4 is devoted to the
study of torsionless-finite algebras, where we show that every
representation-hereditary algebra is torsionless-finite. In Section~5,
we investigate necessary and sufficient conditions for an algebra to
be representation-hereditary. In Section-6, we investigate, from the
perspective of representation-hereditary algebras, the class of
algebras admitting no nonzero morphisms from injective to projective
modules. In Section-7, focusing on tame quasitilted algebras, we
establish necessary conditions for them to be
representation-hereditary and provide sufficient conditions in the
case of tilted algebras.

\section{Conventions and notation}

In this paper, all algebras are assumed to be Artin algebras. For an
algebra $A$, we denote by $\Mod A$ the category of all finitely
generated right $A$-modules, and by $\ind A$, a full subcategory of
$\Mod A$ consisting of exactly one representative from each
isomorphism class of indecomposable modules.  For a category
$\mathcal{C}$, we write $M \in \mathcal{C}$ to indicate that $M$ is an
object in $\mathcal{C}$. We denote by $\add \ \mathcal{C}$ the
subcategory of $\Mod A$ whose objects are finite direct sums of
summands of modules in $\mathcal{C}$. If $M$ is a module, we
abbreviate $\add\{M \}$ as $\add M$. We denote the projective (or
injective) dimension of a module $M$ as $\mathrm{pd} M$ (or
$\mathrm{id} M$, respectively).  We denote by $\mathrm{Gen} M$ (or
$\mathrm{Cogen} M$) the full subcategory of $\mod A$ whose objects are
all modules generated (or cogenerated, respectively) by $M$. We denote
by $\tau = \mathrm{DTr}$ and $\tau^{-1} = \mathrm{Tr D}$ the
Auslander-Reiten translations.

A module $M$ is called a {\it generator} of $\mathrm{mod} A$ if every
projective $A$-module belongs to $\mathrm{add} M$. It is called a {\it
  cogenerator} of $\mathrm{mod} A$ if every injective $A$-module
belongs to $\mathrm{add} M$. Furthermore, it is called a {\it
  generator-cogenerator} of $\mathrm{mod} A$ if it is both a generator
and a cogenerator of $\mathrm{mod} A$. The {\it representation
  dimension} of $A$ is defined as the infimum of the global dimension
of the endomorphism algebra of generator-cogenerator modules in
$\mathrm{mod} A$, that is, $\mathrm{rep.dim} A = \mathrm{inf}
\{\mathrm{gl.dim} \ \mathrm{End}(M) \mid M = A \oplus DA \oplus
M'\}$. A generator-cogenerator $M$ is called an {\it Auslander
  generator} if $\mathrm{gl.dim} \ \mathrm{End}(M) = \mathrm{rep.dim}
A$.

Recall that a morphism $f\colon M \to N$ is said to be {\it right
  minimal} if any morphism $g$ such that $fg = f$ is an
isomorphism. Dually, $f\colon M \to N$ is said to be {\it left
  minimal} if any morphism $g$ such that $gf = f$ is an isomorphism.

Let $\mathcal{X}$ be an additive full subcategory of $\mod A$. For an
$A$-module $M$, a {\it right $\mathcal{X}$-approximation of $M$} is a
morphism $f\colon X \to M$ with $X \in \mathcal{X}$ such that the
sequence of functors $\Hom_A(-, X)|_{\mathcal{X}} \to \Hom_A(-,
M)|_{\mathcal{X}} \to 0$ is exact.

A morphism $f$ is a {\it minimal right $\mathcal{X}$-approximation of
  $M$} if it is both a right $\mathcal{X}$-approximation of $M$ and a
right minimal morphism. Dually, a {\it left
  $\mathcal{X}$-approximation of $M$} is a morphism $f\colon M \to X$
with $X \in \mathcal{X}$ such that the sequence of functors $\Hom_A(X,
-)|_{\mathcal{X}} \to \Hom_A(M, -)|_{\mathcal{X}} \to 0$ is exact. A
morphism $f$ is a {\it minimal left $\mathcal{X}$-approximation of
  $M$} if it is both a left $\mathcal{X}$-approximation of $M$ and a
left minimal morphism.

There are two subcategories of $\mathrm{ind} A$ that deserve special
attention, as they provide, in certain situations, a deeper
understanding of either the module category or, more explicitly, the
algebra itself. To introduce these subcategories, recall that a path in $\mathrm{ind} A$ is a sequence of nonzero morphisms  
$
M = X_{0} \longrightarrow X_{1} \longrightarrow \cdots \longrightarrow X_{t-1} \longrightarrow X_{t} = N$
with $t \geq 1$, where $ X_i \in \mathrm{ind} A $ for all $i$. We then say that $M$ is a \emph{predecessor} of $N$ and that $N$ is a \emph{successor} of $M$.

We denote by ${\mathcal L}$ the full subcategory of
$\mathrm{ind} A$,  defined in \cite{HRS96}, whose objects are the modules $X$ such that, for any
predecessor $Y$ of $X$, the projective dimension $\mathrm{pd}_{A}(Y)$
does not exceed one. This subcategory $\mathcal{L}$ is referred to as
the left part of $\mathrm{mod} A$. Dually, we define the right part
$\mathcal{R}$ of $\mathrm{mod} A$.  In this context, the article
\cite{ACT04} offers valuable insights into the notions we refer to
here.

Let $\mathcal{C}$ be a full additive subcategory of $\mathrm{mod} A$. A module $N \in \mathcal{C}$ is said to be Ext-injective in $\mathcal{C}$ if $\mathrm{Ext}_{A}^1(-, N )|_{\mathcal{C}} = 0$. Dually, the notion of Ext-projective in $\mathcal{C}$ is defined analogously.

\subsection{On Tilting Theory}
We will introduce the necessary background for the topics addressed in
this article. For further details on the subject, see
(\cite[VI-Theorem 3.8]{ASS}).  Recall that a module $T_A \in \md A$ is
called a {\it tilting module} if it satisfies the following
conditions:
\begin{enumerate}
\item $\mathrm{pd}_{A} T \leq 1$;
\item $\Ext^{1}_{A}(T, T) = 0$;
\item The number of pairwise non-isomorphic indecomposable direct
  summands of $T$ equals the rank of the Grothendieck group
  $K_{0}(A)$.
\end{enumerate}

Given a tilting module $T_A$, let $B = \End_A T$. Recall that
$_{B}T_{A}$ induces the torsion pair $(\mathcal{T}(T),
\mathcal{F}(T))$ in $\mod A$ and the torsion pair $(\mathcal{X}(T),
\mathcal{Y}(T))$ in $\mod B$.  Due to the Brenner - Butler Theorem
(\cite{MR675063}), $\mod A$ and $\mod B$ are related as follows: the
functor $\Hom_{A}(T,-): \mod A \rightarrow \mod B$ induces an
equivalence of categories from $\cal T(T)$ to ${\cal Y}(T)$ and the
functor $\Ext_{A}^{1}(T,-): \mod A \rightarrow \mod B$ induces an
equivalence of categories from ${\mathcal{F}}(T)$ to ${\cal X}(T)$.

If $({\cal T}, {\cal F})$ is a torsion pair then $M \in \cal T$ is
Ext-projective in $\cal T$ if and only if $\tau_A M \in \cal F$ and
$N\in \cal F$ is Ext-injective in $\cal F$ if and only if $\tau_A
^{-1} N \in \cal T$ (see \cite[VI.Proposition 1.11]{ASS}).

\subsection{The generator-cogenerator} 
The following result will enable the investigation of aspects related
to the generator-cogenerator in the category of modules and its
connection with certain homological properties.

\begin{theo} \label{CP}
(See \cite{CP04},\cite{EHIS04},\cite{Xi02}) Let $ A $ be an artin
  algebra and $\bar X$ be a generator-cogenerator of $\mathrm{mod} A
  $. The following are equivalent:
\begin{enumerate}
    \item $\rm{gl.dim.} \End \bar{X} \leq 3$.
    \item For all $M \in \ind A \setminus \add \bar X$ there exists an
      exact sequence
    $$\xymatrix{0\ar[r]&X_{1}\ar[r]&X_{0}\ar[r]&M\ar[r]&0}$$ with
      $X_0, X_1 \in \add \bar X$, such that the induced sequence.
    $$\xymatrix{0\ar[r]&\Hom(-,X_{1})\ar[r]&\Hom(-,X_{0})\ar[r]&\Hom(-,M)\ar[r]&0}$$
      is exact in $\add \bar X$.
    \item For all $M \in \ind A \setminus \add \bar X$ there exists an
      exact sequence
    $$\xymatrix{0\ar[r]&M\ar[r]&X^{0}\ar[r]&X^{1}\ar[r]&0}$$ with
      $X^0, X^1 \in \add \bar X$, such that the induced sequence
    $$\xymatrix{0\ar[r]&\Hom(X^{1},-)\ar[r]&\Hom(X^{0}, -
        )\ar[r]&\Hom(M,-)\ar[r]&0}$$ is exact in $\add \bar X$.
\end{enumerate} \theend
\end{theo}

\section{On representation-hereditary algebras: definiton and properties} 

In this section, we introduce the class of representation-hereditary
algebras and present the main result from a homological point of
view. Throughout, we assume that $\Mod A \neq \add(A \oplus DA)$,
meaning that $\Mod A$ contains more modules than just the projective
and injective ones. Otherwise, the global dimension of $\End(A \oplus
DA)$ would be at most 2, since $A$ is representation finite.

\subsection{Representation-hereditary algebras}
\begin{defn}
  An artin algebra $A$, whose $\mathrm{gl.dim} (\End(A \oplus DA)) =
  3$, is called a {\it representation-hereditary algebra}.
\end{defn}

If $A$ is a non-semisimple hereditary algebra, taking $\bar{X} = A
\oplus DA$, it is easy to verify that this algebra satisfies condition
$(2)$ of Theorem \eqref{CP}, and hence $\mathrm{gl.dim} (\End(A\oplus
DA)) \leq 3$. For each $M \in \ind A \setminus \add \bar{X}$, consider
the minimal projective resolution of $M$, $0 \rightarrow
P_{1}\rightarrow P_{0}\rightarrow M \rightarrow 0$, with $P_0, P_1 \in
\add A$. Moreover, if $A$ is assumed to be a representation-infinite
algebra, then $\mathrm{gl.dim} (\End(A\oplus DA)) = 3$
(\cite[Corollary 1.9]{CP04}).

It is also important to mention that in this case, $\mathrm{rep.dim} A
= 3$. This topic is closely related to the present work. Several
highly relevant studies have been published on this subject, and the
ones most pertinent to this work are as follows: \cite{CP04},
\cite{EHIS04}, \cite{Xi02}.

We denote a right $\add (A \oplus DA)$-approximation of $M$ by $f
\colon P_0 \oplus I_0 \to M$ and a left $\add (A \oplus
DA)$-approximation of $M$ by $f \colon M \to P^0 \oplus I^0$, with
$P_0, P^0 \in \add A$ and $I_0, I^0 \in \add DA$. We present a
specific formulation of Theorem \ref{CP} for representation-hereditary
algebras.

\subsection{}
\begin{theo}\label{kernel_projective}
For an artin algebra $A$ such that $\Mod A \ne \add (A \oplus DA)$ the
following are equivalent:
\begin{enumerate}
    \item $A$ is a representation-hereditary algebra.
    \item for all $M \in \ind A \setminus \add (A \oplus DA)$ there
      exists a right $\add (A \oplus DA) $-approximation $f \colon P_0
      \oplus I_0 \to M$ such that $Ker f$ is projective.
    \item for all $M \in \ind A \setminus \add (A \oplus DA)$ the
      kernel of a minimal right $\add (A \oplus DA) $-approximation $f
      \colon P_0 \oplus I_0 \to M$ is projective.
    \item for all $M \in \ind A \setminus \add (A \oplus DA)$ there
      exists a left $\add (A \oplus DA) $-approximation $f \colon M
      \to P^0 \oplus I^0$ such that $Coker f$ is injective.
    \item for all $M \in \ind A \setminus \add (A \oplus DA)$ the
      cokernel of a minimal left $\add (A \oplus DA) $-approximation
      $f \colon M \to P^0 \oplus I^0$ is injective.
\end{enumerate}
\end{theo}

\proof

Suppose that $A$ is a representation-hereditary algebra. By
definition, the global dimension of $\End (A \oplus DA)$ is 3. Then,
by Theorem \ref{CP}, for all $M \in \ind A \setminus \add (A \oplus
DA)$, there exists a right $\add (A \oplus DA)$-approximation $f
\colon P_0 \oplus I_0 \to M$ such that $\Ker f = P_1 \oplus I_1$,
where $P_1$ is a projective module and $I_1$ is an injective module.

Now, consider $\varphi \colon P \oplus I \to M$ to be a minimal right
$\add (A \oplus DA)$-approximation of $M$, and let $K = \Ker
\varphi$. Since both $f$ and $\varphi$ are right $\add (A \oplus
DA)$-approximations, there exist morphisms $\beta_{1}$ and $\beta_{2}$
such that $\varphi = f \beta_{1}$ and $f = \varphi \beta_{2}$. By
taking kernels, we obtain the following commutative diagram with exact
rows:

\[\xymatrix{ 0 \ar[r] & K \ar[r] \ar[d]^{\alpha_{1}} & P \oplus I \ar[d]^{\beta_{1}} \ar[r]^{\varphi} & X \ar@{=}[d] \ar[r] & 0 \\ 0 \ar[r] & P_{1} \oplus I_1 \ar[r] \ar[d]^{\alpha_{2}} & P_{0} \oplus I_{0} \ar[r]^{f} \ar[d]^{\beta_{2}} & X \ar[r] \ar@{=}[d] & 0 \\ 0 \ar[r] & K \ar[r] & P \oplus I \ar[r]^{\varphi} & X \ar[r] & 0 \\ }\]

Since $\varphi$ is right minimal and $\varphi = \varphi
\beta_{2}\beta_{1}$, it follows that $\beta_{2}\beta_{1}$ is an
isomorphism. Consequently, $\alpha_{2}\alpha_{1}$ is also an
isomorphism. As a result, $\alpha_2 \colon P_{1} \oplus I_1 \to K$ is
a split epimorphism, implying that $K$ is a direct summand of $\Ker f
= P_1 \oplus I_1$. Moreover, since $K$ has no injective direct
summand, it follows that $K$ is projective. This proves that
1. implies 2. (and 3.). A similar argument demonstrates that
2. implies 3., while the implication from 3. to 2. is evident. The
implication from 2. to 1. follows from Theorem \ref{CP}, given that
$\Mod A \ne \add (A \oplus DA)$. Therefore, there is an equivalence
among 1., 2., and 3. The remaining statements follow by duality.
\theend

\section{Torsionless-finite algebras}

We define an algebra as {\it torsionless-finite} (see
\cite{zbMATH06154188}) if there exist only finitely many
indecomposable modules, up to isomorphism, that are submodules of
projective modules, or equivalently, if there are finitely many
indecomposable factors of injective modules. It is evident that
torsionless-finite algebras include hereditary ones.

 Ringel proves that $\mathrm{rep.dim} A \leq 3$ provided there are
 only finitely many isomorphism classes of indecomposable $A$-modules
 that are torsionless (i.e., submodules of a projective
 module). However, it is important to note that not all Artin algebras
 with representation dimension at most 3 are torsionless-finite (see
 \cite{zbMATH06154188}).

\begin{theo}\label{torsionless}
  Let $A$ be a representation-hereditary algebra. Then
  \begin{enumerate}
  \item[(a)] If $M$ is a non injective indecomposable module in $\gen
    DA$ then $M = \Omega^{-}(P) $ for some projective indecomposable
    module $P$.
     \item[(b)] If $M$ is a non projective indecomposable module in
       $\cogen A$ then $M = \Omega(I) $ for some injective
       indecomposable module $I$.
\item[(c)] $A$ is torsionless-finite.
\end{enumerate} 
\end{theo}

  \proof

  As $M \in \gen DA$, a right minimal $\add DA$-approximation of $M$
  is an epimorphism $g\colon I \to M$, which is also a minimal right
  $\add (A \oplus DA)$-approximation of $M$. Let $\iota \colon K \to
  I$ denote the kernel of $g$, and let $u \colon K \to I(K)$ be the
  injective envelope of $K$. We have that there exists a split
  epimorphism $t\colon I \to I(K)$ such that $t\iota = u$. Now,
  consider the commutative diagram with exact rows, where $h$ is
  obtained by passing to the cokernels.
  
$$\xymatrix{0 \ar[r] & K \ar@{=}[d] \ar[r]^{\hspace{-0.0cm}\iota}& I
    \ar[d]^{t} \ar[r]^{g}& M \ar[r]\ar[d]^{h}& 0 \\ 0 \ar[r]&
    K\ar[r]_{u} & I(K)\ar[r] & \Omega^-(K) \ar[r] &0}$$

  It follows that $h$ is an epimorphism and that $\mathrm{Ker} h \cong
  \mathrm{Ker} t = I'$, which is injective. Therefore, the short exact
  sequence $0 \to I' \to M \to \Omega^-(K) \to 0$ splits, implying
  that $M \cong I' \oplus \Omega^-(K)$. However, since $M$ is
  indecomposable and non-injective, we must have $M \cong \Omega^-(K)$
  and $I' = 0$.  Furthermore, since the cosyzygy is additive and, by
  Theorem \ref{kernel_projective}, $K$ has no injective direct
  summand, it follows that $K$ is indecomposable and projective. This
  completes the proof of (a). The proof of (b) is analogous, which in
  turn leads to (c) \theend

\section{Some sufficient conditions} \label{sufficient}

  In this section, we present some sufficient conditions for an
  algebra to be representation-hereditary. These conditions arise
  through the study of certain subcategories of the module category.

  For this purpose, we shall define some subcategories.  In the
  remainder of this text, we will denote $\mathcal{P}^{\leq 1} := \{ X
  \in \mathrm{ind}A \mid \mathrm{pd}X \leq 1 \}$ and
  $\mathcal{I}^{\leq 1} := \{ X \in \mathrm{ind}A \mid \mathrm{id}X
  \leq 1 \}$. Recall that the support of a functor $F \colon
  \mathcal{C} \to \mathcal{D}$ is defined as $\mathrm{Supp}F := \{ X
  \in \mathcal{C} \mid F(X) \neq 0 \}$.

  Observe that if $M \in \operatorname{ind}A$ is a non-projective,
  non-injective module, and $\mathrm{Hom}_{A}(DA, M) = 0$, then its
  projective cover provides a minimal right $\mathrm{add}(A \oplus
  DA)$-approximation. Similarly, if $\mathrm{Hom}_{A}(M, A) = 0$, the
  injective envelope of $M$ serves as a minimal left $\mathrm{add}(A
  \oplus DA)$-approximation.


  We now have the following result. The proof is a direct consequence
  of Theorem \ref{kernel_projective}.
  \subsection{}  
 \begin{prop}
 If $A$ is a representation-hereditary algebra then:
 \begin{enumerate}
   \item[(a)] ${\rm indA} \cap {Ker Hom(DA,\_)}\subseteq \mathcal{P}^
     {\le 1}$;
  \item[(b)] ${\rm indA} \cap {Ker Hom(\_, A)}\subseteq \mathcal{I}^
    {\le 1}$.
 \end{enumerate}
 \end{prop}


Recall \ that we assume \ $\Mod A \neq \add (A \oplus DA)$, as the
alternative would imply that $\mbox{$\operatorname{gl.dim} \,
  \operatorname{End}(A \oplus DA) \leq 2$}$.
\subsection{}
\begin{theo} \label{maiprop-1}
Let $A$ be an Artin algebra satisfying the following properties:
\begin{enumerate}
\item[(a.1)] ${\rm ind}A \cap \Ker(\mathrm{Hom}(DA, \_)) \subseteq
  \mathcal{P}^{\leq 1}$;
\item[(a.2)] $\mathrm{Supp}(\mathrm{Hom}(DA, \_)) \setminus {\rm add}A
  \subseteq \mathrm{Gen}(DA)$;
\item[(a.3)] ${\rm ind}A \setminus \mathcal{I}^{\leq 1} \subseteq {\rm
  add}A;$
\end{enumerate}
Then $A$ is a representation-hereditary algebra.
\end{theo}
\proof By Theorem \ref{kernel_projective}, it suffices to prove that
for any indecomposable module in $\ind A \setminus \add (A \oplus
DA)$, its minimal right $\add(A \oplus DA)$-approximation has a
projective kernel. Let $M \in \ind A \setminus \add(A \oplus DA)$ be
an indecomposable module.

First, assume $\mathrm{Hom}(DA, M) = 0$. In this case, the projective
cover of $M$ is itself a minimal right $\add(A \oplus
DA)$-approximation of $M$. By property \textit{(a.1)}, we have $\pd M
= 1$, which implies that the kernel $\Omega(M)$ is projective.

Now suppose $\mathrm{Hom}(DA, M) \neq 0$. Then $M \in
\mathrm{SuppHom}(DA, \_)$, and by property \textit{(a.2)}, it follows
that $M \in \gen(DA)$. By the proof of Theorem \ref{torsionless},
there exists an indecomposable module $K$ such that $M \cong
\Omega^-(K)$. Since $M$ is non-injective, it follows that $\id K \geq
2$. By property \textit{(a.3)}, this implies that $K$ is
projective. Using again the proof of Theorem \ref{torsionless}, $K$ is
the kernel of a minimal right $\add(DA)$-approximation of $M$, which
in this case is also a minimal right $\add(A \oplus DA)$-approximation
of $M$. This completes the proof that $A$ is a
representation-hereditary algebra.  \theend

We also have the dual version of this theorem, whose proof is similar.
\subsection{}
\begin{theo} \label{suficiente}
Let $A$ be an Artin algebra satisfying the following properties:
\begin{enumerate}
\item[(b.1$ $)] $ {\rm ind}A \cap \Ker(\mathrm{Hom}(\_, A)) \subseteq
  \mathcal{I}^{\leq 1}$;
\item[(b.2)] $\mathrm{ind}A \setminus \mathcal{P}^{\leq 1} \subseteq
  \mathrm{add}DA;$
\item[(b.3)] $\mathrm{Supp}(\mathrm{Hom}(\_, A)) \setminus
  \mathrm{add}DA \subseteq \mathrm{Cogen}A$;
\end{enumerate}
Then $A$ is a representation-hereditary algebra.
\end{theo}

As a consequence of the preceding results, we provide sufficient conditions on the subcategories $\mathcal{L}$ and $\mathcal{R}$ to classify algebras as representation-hereditary.



\subsection{}
\begin{cor}
Let $A$ be an artin algebra.
\begin{enumerate}
    \item[(a)] If $\ind A \setminus {\mathcal L} \subseteq \{ I \mid I
      \mbox{ is simple injective} \}$ and $\ind A \setminus {\mathcal
        R} \subseteq \add (A \oplus DA)$ \\ then $A$ is a
      representation-hereditary algebra.

    \item[(b)] If $\mbox{ind}A\setminus\mathcal{R}\subseteq
      \{P\;|\;P\mbox{ is a simple projective } \}$ and
      $\mbox{ind}A\setminus\mathcal{L}\subseteq \mbox{add} (A \oplus
      DA)$ then $A$ is a representation-hereditary algebra.
\end{enumerate}

\end{cor}
\proof According to the hypothesis $\mbox{ind}A\setminus
\mathcal{L}\subseteq\{I\;|\;I\mbox{ is injective simple}\}\subseteq
\mbox{add}DA$, then $\mbox{ind}A\setminus \mbox{add}DA\subseteq
\mathcal{L}$ and then $({\rm indA}\setminus {\rm add}DA)\setminus {\rm
  SuppHom(DA,\_)}\subseteq\mathcal{L} \subseteq \mathcal{P}^ {\le
  1}$. Furthermore $\mbox{ind}{A}\setminus \mathcal{I}^{\le
  1}\subseteq \mathrm{ind}A\setminus \mathcal{R} \subseteq \add (A
\oplus DA),$ then $\mbox{ind}{A}\setminus \mathcal{I}^{\le 1}\subseteq
\mathrm{add}A.$ These proves the items {\it (a.1)} and {\it (a.3)} in
Theorem \ref{maiprop-1}.  In order to prove {\it (a.2)}, let $X\in
\mbox{ind}A\setminus \add (A \oplus DA)$ such that
$\mbox{Hom}(DA,X)\ne 0.$

As $X$ is not a simple injective module by hypothesis $X\in
\mathcal{L}$, then there exists a non zero sectional path of
irreducible maps from some injective $I$ to $X$ in ind$A,$ (see
\cite{ACT04}, 3.1, 3.2) says

$$\xymatrix{I\ar[r]^{f}&X_{0}\ar[r]&\ldots\ar[r]&X_{n-1}\ar[r]&X_{n}=X}$$
We will suppose $n$ minimal in the sense that for each $i\in
\{0,1,2,\ldots,n\}$ the module $X_{i}$ is not an injective module. It
is known that for each $X_{i},$ Hom$(DA,X_{i})\ne 0,$ that is,
$\tau^{-}X_{i}\notin \mathcal{L}$, then $\tau^{-}X_{i}$ is a simple
injective $A$-module.  Since $X_{0}\rightarrow X_{1}$ is irreducible,
there is a irreducible morphism $X_{1}\rightarrow\tau^{-}X_{0}$ which
$\tau^{-}X_{0}$ is a simple injective. Then $X_{1}$ must be a
injective module in contradiction with the minimality of $n.$ Thus
$X=X_{0}$ and $f$ is a irreducible map. Hence $f$ is an epimorphism
and consequently $X\in \gen(DA).$ The result of item (a) follows by
Theorem \ref{maiprop-1}. The proof of (b) is dual.  \theend

\subsection{}
\begin{cor}
If $A$ is an artin algebra such that $$\mathrm{ind}A\setminus
(\mathcal{L}\cap \mathcal{R})\subseteq \{\mbox{simple projective or
  simple injective}\},$$ then $A$ is a representation-hereditary
algebra.
\end{cor}
\proof Clearly, $\mathrm{ind}A \setminus \mathcal{L} \subseteq
\mathrm{ind}A \setminus (\mathcal{L} \cap \mathcal{R})$.  Observe that
any simple projective module lies in $\mathcal{L}$, and thus
$\mathrm{ind}A \setminus \mathcal{L} \subseteq \{\text{simple
  injective modules}\}$.  Consequently, the result follows from the
previous corollary.

\theend

\section{Algebras with no morphism from injective to projective} \label{no injective to projective}

A natural question that arises regarding representation-hereditary 
algebras concerns their global dimension. The following example
demonstrates that there is no upper bound for the global dimension of
these algebras. The remainder of this section will be devoted to
presenting some results concerning representation-hereditary algebras
in the case where $\mathrm{Hom}_{A}(DA, A) = 0$.

\subsection{}

The  algebra $A$ given by the following quiver  
$$  
\xymatrix{  
1 \ar@(ul,dl)[]_{\alpha} \ar[r]^{\beta} & 2  
}  
$$  
and bound by the relations $\{ \alpha^{2}, \alpha \beta \}$ is of finite representation type and infinite global dimension.

Computing the Auslander-Reiten quiver of $A$ yields 
{\scriptsize
\begin{center}
\begin{tikzcd}[column sep=small, row sep=small]
& P_2 \ar{dr} &  & I_1 \ar{dr}  &  &  &   \\
   &  & P_1 \ar{ur} \ar{dr}  & &  S_1 & &  \\
& S_1 \ar{ur} &  & I_2 \ar{ur}  &   & & .
\end{tikzcd}
\end{center}
}
The only module that is neither projective nor injective is the simple module $S_1$, and the corresponding sequences given by the left and right $\mathrm{add}(A \oplus DA)$-approximations of this module are clearly depicted in the quiver. According to Theorem \eqref{kernel_projective}, this algebra is representation-hereditary.

\subsection{}
\begin{prop} \label{morfismodeDAtauxdeterminaparax}
Let $A$ be a representation-hereditary algebra.
\begin{itemize}
    \item[a)] \label{orbita} If $X \in \mathrm{ind}A$ is such that
      $\Hom(DA, \tau X) \neq 0$, then $\Hom(DA, X) \neq 0$.
    \item[b)] \label{orbita2} If $X \in \mathrm{ind}A$ is such that
      $\Hom(\tau^- X, A) \neq 0$, then $\Hom(X, A) \neq 0$.
\end{itemize}
\end{prop} 
\proof We prove only \textit{(a)} because \textit{(b)} is its
dual. Let $X \in \ind A$ be such that $\Hom(DA, \tau X) \neq 0$. In
particular, $X$ is non-projective and satisfies $\pd X \geq 2$ (see
\cite{ASS}, (IV.2.7)).

Suppose, by contradiction, that $\Hom(DA, X) = 0$. Then $X$ is
non-injective, and the projective cover of $X$ is a minimal right
$\add(A \oplus DA)$-approximation of $X$. Since $A$ is a
representation-hereditary algebra, by Theorem~\ref{kernel_projective},
the kernel $\Omega(X)$ is projective.

This leads to a contradiction because $\pd X \geq 2$. Therefore, we
conclude that $\Hom(DA, X) \neq 0$.  \qed

\subsection{} 
\begin{cor} \label{dimensaoprojetivainjetiva}
Let $A$ be a representation-hereditary algebra.
\begin{itemize}
    \item[a)] If $X$ is an indecomposable non-injective module with $\pd X \geq 2$, then $\pd  \tau^{-} X \geq 2$.
    \item[b)] If $X$ is an indecomposable non-projective module with $\id X \geq 2$, then $\id      \tau X \geq 2$.
\end{itemize}
\end{cor}
\qed

This last result implies that if some non injective indecomposable
module has projective dimension at least two then any module in the
same $\tau^{-}$-orbit has also projective dimension at least
two. Dually, if some non projective module has injective dimension at
least two then any module in the same $\tau$-orbit has also injective
dimension at least two.

\subsection{}
Following
\cite{MR1648603}, an algebra is called \emph{quasitilted} if $A =
\mathrm{End}, T$, where $T$ is a tilting object in a
$\mathrm{Hom}$-finite abelian $k$-hereditary category. In
\cite{HRS96}, it is shown that this definition is equivalent to
requiring that the algebra has global dimension at most two and that,
for any indecomposable module $M$, either $\pd M \leq 1$ or $\id M
\leq 1$.

\begin{theo} \label{quasitilted}
  Let $A$ be a representation-hereditary algebra. If $\Hom(DA, A) = 0$, then one of the following conditions hold:
  \begin{enumerate}
  \item[(a)] $A$ is a quasitilted algebra or
   \item[(b)] There exist an indecomposible module $X$ such that  $\mathrm{pd}\, X = \mathrm{id}\, X = 2$, then $  \mathrm{pd}\, \tau^{-1} X = \mathrm{id}\, \tau X = 2.$
   \end{enumerate}
  \end{theo}
\proof Let $P$ be a projective module and $X \in \ind A$ be a non
projective summand of rad$P$. As $\Hom(DA,P) = 0$ then $\Hom(DA, X) =
0$ and so $\Hom(DA, \tau X) = 0$, because Proposition
\ref{orbita}. This implies that $\pd X \leq 1$ and so $\pd
(\mathrm{rad}P) \leq 1$ for any projective $P$.  By the short exact
sequence
$$0 \to \mathrm{rad} P \to P \to \mathrm{top} P \to 0$$ we conclude
$\pd S \leq 2$ for any simple module $S$ and so gl.dim$A \leq 2$. 

  Now, let $X \in \text{ind}A$ be such that $\text{pd}X = \text{id}X =
  2$. Since $\text{pd}X = 2$, it follows from Proposition
  \eqref{morfismodeDAtauxdeterminaparax} that $\text{Hom}(DA, X) \neq
  0$. Similarly, since $\text{id}X = 2$, we conclude that
  $\text{Hom}(X, A) \neq 0$. \qed

\begin{prop} \label{quasitilted}
  Let $A$ be a representation-hereditary algebra. If $\Hom(DA, A) = 0$, then for any indecomposable module $X$ we have:
  \begin{enumerate}
    \item[(a)] $ X \notin \operatorname{gen} DA $, and the morphism
    $( f\, p) \colon I \oplus P_0(C) \to X$  
    is a minimal right $\ensuremath{\operatorname{add}(A \oplus DA)}$-approximation
    of $ X $, where:
    \begin{enumerate}
      \item[(a$_{1}$)] $ f \colon I \to X $ is a minimal right
        $\ensuremath{\operatorname{add}(DA)}$-approximation of $ X $;
      \item[(a$_{2}$)] $ c \colon X \to C $ is the cokernel of $ f $;
      \item[(a$_{3}$)] $ \pi \colon P_0(C) \to C $ is the projective
        cover of $ C $;
      \item[($a_{4}$)] $ p \colon P_0(C) \to X $ satisfies $ c \, p =
        \pi $.
    \end{enumerate}
    \item[(b)] $ X \notin \operatorname{cogen} A $, and the morphism
    $ \begin{bmatrix}
    f \\ p
    \end{bmatrix}
    \colon X \to P \oplus I_0(K)$    
    is a minimal left $\ensuremath{\operatorname{add}(A \oplus DA)}$-approximation
    of $ X $, where:
    \begin{enumerate}
      \item[(b$_{1}$)] $ f \colon X \to P $ is a minimal left
        $\ensuremath{\operatorname{add}(A)}$-approximation of $ X $;
      \item[(b$_{2}$)] $ k \colon K \to X $ is the kernel of $ f $;
      \item[(b$_{3}$)] $ \iota \colon K \to I_{0}(K) $ is the injective
        envelope of $ K $;
      \item[($b_{4}$)] $ i \colon X \to I_{0}(K) $ satisfies $ i \, k =
        \iota $.
    \end{enumerate}
  \end{enumerate}
\end{prop}
\proof Let $f: I \to X$ be a minimal right $\text{add}DA$-approximation and
$g: X \to P$ a minimal left $\text{add}A$-approximation of $X$. By
hypothesis, $gf = 0$, implying that $f$ is not an epimorphism, and $g$
is not a monomorphism. Thus, $X \notin \gen DA$ and $X \notin \cogen
A$.

Observe that if $P$ is projective and $\left(f g\right)\colon I \oplus
P \to X$ is an epimorphism, then $\left(f g\right)$ is a right
$\text{add}\left(A \oplus DA\right)$-approximation of $X$. This
follows because any morphism from an injective object to $X$ factors
through $f$, and any morphism from a projective object to $X$ factors
through any epimorphism.

Let $c \colon X \to C$ be the cokernel of $f$, and let $\pi \colon
P_0(C) \to C$ denote the projective cover of $C = \mathrm{Coker}
f$. Then there exists a morphism $p \colon P_0(C) \to X$ such that $c
\, p = \pi$. The following commutative diagram with exact rows
illustrates this, where $f'\colon I \to \mathrm{Im} f$ is the
corestriction of $f$, and $j\colon \mathrm{Im} f \to X$ is the
inclusion:

\[
\begin{small}
\begin{array}{c}
\xymatrix@C=1cm@R=0.8cm{ 0 \ar[r] & I
  \ar[r]^{\hspace{-0.5pc}\raisebox{2pt}{\scriptsize $\begin{pmatrix} 1
        \\ 0 \end{pmatrix}$}} \ar[d]_{\scriptsize f'} & I \oplus
  P_0(C) \ar[r]^{\hspace{0.5pc}\raisebox{2pt}{\scriptsize
      $\begin{pmatrix} 0 & 1 \end{pmatrix}$}}
  \ar[d]^{\raisebox{2pt}{\scriptsize $\begin{pmatrix} f &
        p \end{pmatrix}$}} & P_0(C) \ar[r] \ar[d]_{\scriptsize \pi} &
  0 \\ 0 \ar[r] & \operatorname{Im} f \ar[r]_{\scriptsize j} \ar[d] &
  X \ar[r]_{\scriptsize c} & C \ar[r] \ar[d] & 0 \\ & 0 & & 0 }
\end{array}
\end{small}
\]

As $f'$ and $\pi$ are epimorphisms, $\left(f\, p\right)$ is also an
epimorphism and thus a right add$\left(A \oplus
DA\right)$-approximation of $X$. As $\Hom(DA,A) = 0$, we prove that
the right add$\left(A \oplus DA\right)$-approximation $\left(f
p\right) \colon I \oplus P_0(C) \to X$ is minimal.

If $\Hom(I,P_0(C)) = 0$ then the minimal right add$\left(A \oplus
DA\right)$-approximation of $X$ is $(f , g)\colon I \oplus P \to X \to
0$ for some projective $P$ such that $P_0(C) = P \oplus P'$, $g =
p\mid_P$ and $p\mid_{P'} = 0$. Now, consider the commutative diagram
with exact lines where $h = cg$ and $f = jf'$ is the canonical
decomposition of $f$:

\[
\begin{small}
\begin{array}{c}
\xymatrix@C=1cm@R=0.8cm{ 0 \ar[r] & I
  \ar[r]^{\raisebox{2pt}{\scriptsize $\begin{pmatrix} 1
        \\ 0 \end{pmatrix}$}} \ar[d]_{\scriptsize f'} & I \oplus P
  \ar[r]^{\raisebox{2pt}{\scriptsize $\begin{pmatrix} 0 &
        1 \end{pmatrix}$}} \ar[d]^{\raisebox{2pt}{\scriptsize
      $\begin{pmatrix} f & g \end{pmatrix}$}} & P \ar[r]
  \ar[d]_{\scriptsize h} & 0 \\ 0 \ar[r] & \operatorname{Im} f
  \ar[r]_{\scriptsize j} \ar[d] & X \ar[r]_{\scriptsize c} \ar[d] & C
  \ar[r] & 0 \\ & 0 & 0 & }
\end{array}
\end{small}
\]

Then $h\colon P \to C$ is an epimorphism and so it factors through the
projective cover $\pi \colon P_0(C) \to C$. By minimality of $\pi$ we
have that $P_0(C)$ is a direct summand of $P$. Therefore $P = P_0(C)$
and $g = p$. \qed

\section{On representation-hereditary quasi-tilted algebras}

A general study of tame quasitilted algebras can be found in
\cite{MR1622799}. A particular subclass of quasitilted algebras is
that of tilted algebras. Recall that a \emph{tilted} algebra is an
algebra of the form $A = \End_H(T_H)$, where $T_H$ is a tilting
module, meaning that $\Ext_H(T_H, T_H) = 0$ and the number of pairwise
non-isomorphic indecomposable direct summands of $T_H$ equals the
number of vertices of the quiver $Q_H$, with $H = kQ$ a hereditary
algebra (see \cite{ASS}).  If $Q$ is a Dynkin quiver, then $A$ is of
finite representation type. If $Q$ is Euclidean, we say that $A$ is of
Euclidean type; otherwise, $A$ is said to be of wild type.

If $H = kQ$, where $Q$ is an acyclic quiver whose underlying graph is
Euclidean, then the tilted algebra $A$ is of finite representation
type if and only if $T$ has both a postprojective and a preinjective
summand (see \cite[XVII-Lemma 3.3]{MR2382332}).  If $H$ is
representation-infinite and $T_H$ is a postprojective (or
preinjective) module, then $A$ is called a \emph{concealed} algebra
(see \cite[Theorem 4.5, Theorem 4.7]{ASS}).  In the particular case
where $H$ is of Euclidean type, $A$ is referred to as a tame concealed
algebra.  If $A$ is not concealed, then $T$ necessarily has a regular.

It is worth noting that a component \( \Gamma \) of the Auslander-Reiten quiver is called semiregular if \( \Gamma \) does not contain both a projective and an injective module at the same time. If \( \Gamma \) is a component of the Auslander-Reiten quiver of a quasi-tilted algebra containing an oriented cycle, then \( \Gamma \) is a semiregular tube (\cite[Theorem A]{MR1372358}).

Recall that a semiregular tube is either regular, meaning a stable tube, or it is obtained from a regular tube by a finite sequence of ray (or coray) insertions (see \cite[XV.2]{MR2382332}. Another important result states that if a quasi-tilted algebra is not tilted, then every component of its Auslander-Reiten quiver is semiregular. The structure of the components of the Auslander-Reiten quiver of quasi-tilted algebras was studied in \cite{MR1372358}. The structure of the Auslander-Reiten quiver of tame tilted algebras is presented in \cite[Theorem 3.5, 3.6]{MR2382332} and \cite{MR989917}. Succinctly, these results establish that a representation-infinite tilted algebra of Euclidean type is either a domestic tubular extension or a domestic tubular coextension of a concealed algebra of Euclidean type, there are a branch enlargement of a tame concealed algebra.  In \cite{MR989917} the Auslander-Reiten component of tame tilted of wild type are described.

\subsection{}
\begin{theo} \label{concealed}
  Let $A$ be a representation-hereditary algebra.  
  If $A$ is a tame quasi-tilted algebra, then $A$ is tilted of one of the following types:
  \begin{enumerate}
   \item[(a)] $A$ is tilted of finite representation type.
   \item[(b)] $A$ is tame concealed.
   \item[(c)] $A$ is tame tilted of wild type, without inserted or coinserted tubes, and the slice of $A$ lies in the connecting component, which is neither a postprojective nor a preinjective component of $A$.
  \end{enumerate}
\end{theo}

 \proof  If $A$ is representation-finite, then $A$ is tilted (\cite[II-Corollary 3.6]{HRS96}). If $A$ is of infinite representation type, we first prove that, in general, representation-hereditary algebras do not admit standard tubes obtained by ray or coray insertions and, consequently, do not admit standard tubes obtained by a finite sequence of ray or coray insertions, which is equivalent to a branch enlargement (\cite[XX-3.23, XV.2, XV-Theorem 3.9]{MR2382332}; \cite[4.7]{Rin84}).

As a consequence, the Auslander-Reiten components of representation-hereditary quasi-tilted algebras that contain an oriented cycle are exclusively regular components. Since we are in the tame case, these components must necessarily be regular tubes.

In the proof of  \cite[Theorem 2.5]{MR1214905}, we have that the structure of a coray tube (dually, for a ray tube) can be described by an infinite sectional path of the form  
$$ \cdots \rightarrow \tau^{2r} X_{1} \rightarrow \tau^{r} X_{s} \rightarrow \cdots \rightarrow \tau^{r} X_{1} \rightarrow X_{s} \rightarrow \cdots \rightarrow X_{1}.$$

From the proof, we have that the modules $X_{i}$  are left stable and belong to the 
$\tau$-orbit of an injective module. Moreover, since these modules are left stable, part of the structure of this component is determined by considering the Auslander-Reiten translate of these modules along the sectional path. Thus, up to a finite number of modules, the figure below represents the general structure of a coray tube.

{\tiny
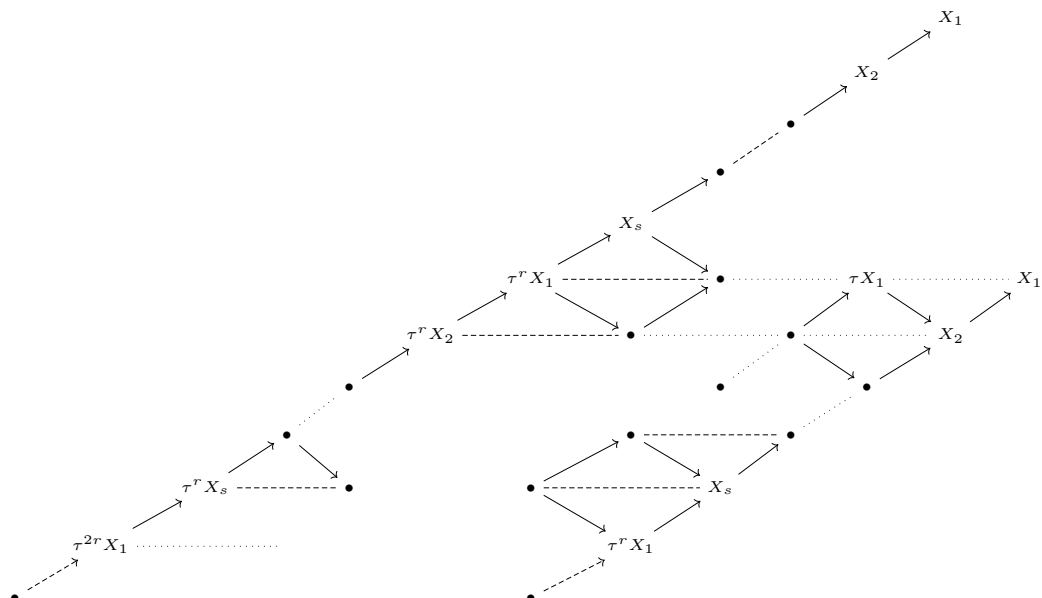
\begin{tikzcd}
	&&&&&&&&&&& {X_{1}} \\
	&&&&&&&&&& {X_{2}} \\
	&&&&&&&&& \bullet \\
	&&&&&&&& \bullet \\
	&&&&&&& {X_{s}} \\
	&&&&&& {\tau^{r}X_{1}} && \bullet && {\tau X_{1}} && {X_{1}} \\
	&&&&& {\tau^{r}X_{2} } && \bullet && \bullet && {X_{2}} \\
	&&&& \bullet &&&& \bullet && \bullet \\
	&&& \bullet &&&& \bullet && \bullet \\
	&& {\tau^{r} X_{s}} && \bullet && \bullet && {X_{s}} \\
	& {\tau^{2r} X_{1}} && {} &&&& {\tau^{r} X_{1}} \\
	\bullet &&&&&& \bullet
	\arrow[from=2-11, to=1-12]
	\arrow[from=3-10, to=2-11]
	\arrow[dashed, no head, from=4-9, to=3-10]
	\arrow[from=5-8, to=4-9]
	\arrow[from=5-8, to=6-9]
	\arrow[from=6-7, to=5-8]
	\arrow[dashed, no head, from=6-7, to=6-9]
	\arrow[from=6-7, to=7-8]
	\arrow[dotted, no head, from=6-9, to=6-11]
	\arrow[dotted, no head, from=6-11, to=6-13]
	\arrow[from=6-11, to=7-12]
	\arrow[from=7-6, to=6-7]
	\arrow[dashed, no head, from=7-6, to=7-8]
	\arrow[from=7-8, to=6-9]
	\arrow[dotted, no head, from=7-8, to=7-10]
	\arrow[from=7-10, to=6-11]
	\arrow[dotted, no head, from=7-10, to=7-12]
	\arrow[from=7-10, to=8-11]
	\arrow[from=7-12, to=6-13]
	\arrow[from=8-5, to=7-6]
	\arrow[dotted, no head, from=8-9, to=7-10]
	\arrow[from=8-11, to=7-12]
    	\arrow[dotted, no head, from=9-4, to=8-5]
	\arrow[from=9-4, to=10-5]
	\arrow[dashed, no head, from=9-8, to=9-10]
	\arrow[from=9-8, to=10-9]
	\arrow[dotted, no head, from=9-10, to=8-11]
	\arrow[from=10-3, to=9-4]
	\arrow[dashed, no head, from=10-3, to=10-5]
	\arrow[from=10-7, to=9-8]
	\arrow[dashed, no head, from=10-7, to=10-9]
	\arrow[from=10-7, to=11-8]
	\arrow[from=10-9, to=9-10]
	\arrow[from=11-2, to=10-3]
	\arrow[dotted, no head, from=11-2, to=11-4]
	\arrow[from=11-8, to=10-9]
	\arrow[dashed, from=12-1, to=11-2]
	\arrow[dashed, from=12-7, to=11-8]
\end{tikzcd}
}

By \cite{MR1414820}, the family of coray tubes is standard. Given the structure of the quiver of a coray tube, there exists an indecomposable injective module $I$ within the tube and a nonzero morphism from $I$ to $\tau X$ for some indecomposable module $X$
 in the tube. By Proposition \eqref{morfismodeDAtauxdeterminaparax}, it follows that there is also a nonzero morphism from $I$ to $X$. 
 However, since coray tubes are standard, the structure of the coray tube ensures that no morphisms exist from  $I$  to the indecomposable module $X$ within the same tube. The remaining injective modules belong to the preinjective component, and it is well known that there are no morphisms from preinjective component to a coray tube.

 If $A$ is a connected tame quasi-tilted algebra, then, by \cite[Theorem A]{MR1622799} and \cite[XX-Theorem 3.26]{MR2382332}, $A$ is either a tame tilted algebra or a tame semiregular branch enlargement of a tame concealed algebra. As we have proven above, since we do not have branch enlargement, then in the second case we have a tame concealed algebra (see also \cite[XVII-Theorem 3.5, 3.6]{MR2382332}).  Let us analyze the case where $A$ is tame tilted.

 Since $A$ is tame tilted and of infinite representation type, it is either of Euclidean type or of wild type. If $A$ is of Euclidean type, then it is a  tubular extension (or coextension) of a concealed algebra of Euclidean type (\cite[XVII-Theorem 5.1]{MR2382332}).
However, as already proven, there are no standard tubes obtained by ray or coray insertions.

For the wild case, we will use the following fact: an algebra is tilted if and only if it admits a slice (\cite[Theorem 5.6-VIII]{ASS}). For the wild case, we can state that the slice lies in the connecting component, which is neither a postprojective component nor a preinjective component. In fact, if the slice were contained in a preinjective component, then by the Brenner-Butler Theorem, $A$ would be wild. Similarly, if it were contained in a postprojective component.
\qed

\subsection{On representation-hereditary  tilted algebras}

In this section, we present sufficient conditions for a tilted algebra
to be representation-hereditary. To this end, we introduce some
additional notation specific to tilted algebras, along with the
foundational background required for a proper understanding of the
topic. We adopt a notation that explicitly indexes each algebra. For
instance, given an algebra $A$, we denote by $\mathcal{L}_A$ the left
part of $\Mod A$, and by $P_A(i)$ (respectively, $I_A(i)$) the
projective (respectively, injective) indecomposable module associated
with the vertex $i$ in the ordinary quiver $Q_A$ of $A$.

Given a tilting module $T_H$ and $A = \mathrm{End}(T)$, we denote by
$(\mathcal{T}(T), \mathcal{F}(T))$ and $(\mathcal{X}(T),
\mathcal{Y}(T))$ the classical torsion pairs induced by $T_H$ in $\Mod
H$ and $\Mod A$, respectively. We also denote by $t$ the torsion
radical of $(\mathcal{T}(T), \mathcal{F}(T))$. If $H$ is hereditary,
it is known that the pair $(\mathcal{X}(T), \mathcal{Y}(T))$ splits,
and that $\mathcal{Y}(T) \subseteq \mathcal{L}$ (see \cite[Lemma
  5.1]{MR675063}).

By \cite[Theorem 3.5]{ASS}, the direct sum of the $A$-modules $E =
\Hom_H(T, DH)$ forms a section lying in an acyclic component of the
Auslander-Reiten quiver of $A$. Moreover, all predecessors of $E$ in
this component belong to $\mathcal{Y}(T)$, and any proper successor of
$E$ in this component belongs to $\mathcal{X}(T)$. According to
\cite[VI-Lemma 4.9]{ASS}, we have that $ \tau^{-} \Hom_{H}(T, I(i))
\simeq \Ext^{1}_{H}(T, P(i))$. Since the successors of $\Hom_{H}(T,
I(i))$ in the Auslander-Reiten component belong to $\mathcal{X}(T)$,
it follows that $\Ext^{1}_{H}(T, P(i))$ is in $\mathcal{X}(T)$. On the
other hand, $\Hom_{H}(T, I(i)) \in \mathcal{Y}(T)$, so we conclude
that $\Hom_{H}(T, I(i))$ is $\Ext$-injective in $\mathcal{Y}(T)$. In
the next result, we will prove sufficient conditions for these to be
all the $\Ext$-injective objects in $\mathcal{L}_A$.

\subsection{}
\begin{lem}\label{slice_L}
Let $A$ be a tilted algebra, $H$ a hereditary algebra, and $T_H$ a tilting module such that $A = \End T$.  

The $A$-module $E = \Hom_H (T, DH)$ is the direct sum of all $\Ext$-injective modules in the left part $\La_A$ if and only if, for every sink $i$ in the ordinary quiver $Q_H$, the projective (simple) module $P_H(i)$ belongs to $\add T$.
\end{lem}
\proof According to \cite[Theorem 3.1]{ACT04}, the $\Ext$-injective modules in $\mathcal{L}_{A}$ are the modules $X \in \mathcal{L}_{A}$ for which there exists a path from an injective module to $X$, or the modules $X$ for which there are no injective modules $I$ with a path to $X$, but there exists a projective module $P \in \ind A \setminus \mathcal{L}_{A}$ and a sectional path from $P$ to $\tau^{-} X$.  

However, all projective modules belong to $\mathcal{L}_{A}$ when dealing with a tilted algebra (\cite[II-Theorem 1.14]{HRS96}). Thus, since $A$ is tilted, the $\Ext$-injective modules in $\La_A$ are those modules $X \in \La_A$ such that there is a path from an injective module $I$ to $X$.

Suppose that $E = \Hom_H (T, DH)$ is the direct sum of the Ext-injective modules of the left part $\La_A$. Let $i$ be a vertex in $Q_H$, and let $I_H(i)$ be the corresponding indecomposable injective module in $\mod H$. Then, there exists an injective module $I$ in $\mod H$ and a path from $I$ to $\Hom_H(T, I_H(i))$ in $\mathcal{L}_A$.  

If $i$ is a sink, let us show that $P_H(i) \in \add T$. Indeed, since there is a path from $I$ to $\Hom_H(T, I_H(i))$ and $I$ is injective in $\mathcal{L}_A$, it follows that it is $\Ext$-injective. Hence, $I = \Hom_H(T, I_H(j))$, so there is a path from $I_H(j)$ to $I_H(i)$, and thus a path from $i$ to $j$ in the ordinary quiver $Q_H$. If $i$ is a sink, then $I = \Hom_H(T, I_H(i))$ is injective. However, we know from \cite[VI-Lemma 4.9]{ASS} that the projective module $P_H(i) \in \add T$ if and only if $\Hom_H(T, I_H(i))$ is an injective module in $\La_A$.

On the other hand, we know that $E = \Hom_H (T, DH)$ is the $\Ext$-injective module in the torsion-free class $\mathcal{Y}(T)$, and we have that $\mathcal{Y}(T) \subseteq \La_A$. In particular, $E$ belongs to $\add \La_A$. Let $X_A = \Hom_H(T, I_H(i))$ for some $i \in (Q_H)_0$. If $i$ is a sink, then by hypothesis, $P_H(i) \in \add T$, so $X_A$ is an injective module in $\La_A$ and thus $\Ext$-injective in $\La_A$. If $i$ is not a sink, then there is a path in $Q_H$ (which is a tree) from $i$ to some sink $j$. Consequently, there is a path from $I_H(j)$ to $I_H(i)$ in $\mathcal{T}(T) \subseteq \Mod H$. Since $j$ is a sink, again by hypothesis, the module $\bar{I}_A = \Hom_H(T, I_H(j))$ is injective, and there is a path from $\bar{I}_A$ to $X_A$. This proves that $X_A$ is $\Ext$-injective in $\La_A$.

\theend

\subsection{}
\begin{lem}\label{cogen_tau}
Let $A$ be a tilted algebra, $H$ be a hereditary algebra and $T_H$ be a tilting module such that $A = \End T$ and $Hom_H(T, DH)$ is the $\Ext$-injective module in $\La_A$.
Then $\ind A \setminus \La_A \subseteq \add DA$ if and only if $\cogen (\tau_H T) = \add (\tau_H T)$.
\end{lem}
\proof
If $T$ is a tilting module in $\mod H$, then $\mathcal{F}(T) = \cogen (\tau_H T)$ (see \cite[VI-Theorem 2.5]{ASS}). The modules in $\ind A \setminus \La_A$ are precisely the modules $\Ext^1_H (T, N)$ with $N_H \in \cogen (\tau_H T)$. However, $\Ext^1_H (T, N)$ is injective if and only if $N_H \in \add (\tau_H T)$ (see \cite[VI-Proposition 5.8]{ASS}).
\theend

\subsection{}
\begin{lem}  
Let $A$ be a tilted algebra, $H$ a hereditary algebra, and $T_H$ a tilting module such that $A = \End T$ and $\Hom_H(T, DH)$ is the $\Ext$-injective module in $\La_A$.  Define $P$ as the direct sum of all projective modules that do not belong to $\add T$, and let $P'$ be the direct sum of the projective modules in $\add T$.  

Then, the inclusion $\ind A \setminus \La_A \subseteq \add DA$ holds if and only if  
\[
T = \tau_H^- \left(\frac{P}{tP}\right) \oplus P'.
\]
\end{lem}
\proof
Suppose $\ind A \setminus \La_A \subseteq \add DA$.
If $P_H(i)$ is a projective $H$-module not lying in $\add T$, then $\Hom_H (T, I_H(i)) \in \La _A$ is not an injective $A$-module. Consequently, by hypothesis, $\tau_A^- \Hom_H (T, I_H(i)) = \Ext^1 _H (T,P_H(i))$ is not in $\La_A$ and is therefore an injective $A$-module.  

Furthermore, we have $\Ext^1_H (T,P_H(i)) = \Ext^1 _H (T,\tau_H T')$ for some direct summand $T_H'$ of $T_H$. It follows from the canonical exact sequence of $P_H(i)$, relative to the torsion pair $(\mathcal{T}(T), \mathcal{F}(T))$, that  
\[
\Ext^1 _H \left( T,P_H(i) \right) = \Ext^1 _H \left( T,\frac{P_H(i)}{tP_H(i)} \right).
\]  

By the Brenner-Butler theorem, it follows that  
\[
\frac{P_H(i)}{tP_H(i)} = \tau_H T'.
\]  
Therefore, we conclude that  
\[
T_H' = \tau_H ^-\left(\frac{P_H(i)}{tP_H(i)}\right).
\]

Reciprocally, suppose that $T = \tau_H^- \left(\frac{P}{tP}\right) \oplus P'$. Any $A$-module in $\ind A \setminus \La_A$ has an injective successor. However, the Ext-projective modules in $\ind A \setminus \La_A$ are precisely those of the form $\Ext^1 _H (T,P_H(i))$ with $P_H(i) \notin \add T$.  

Since  
\[
\Ext^1 _H (T,P) = \Ext^1 _H \left(T,\frac{P}{tP}\right) = \Ext^1 _H (T,\tau_H T)
\]  
is injective, the remaining injective $A$-modules lie in $\La_A$. Consequently, there are no other modules in $\ind A \setminus \La_A$.

\theend

\subsection{}
\begin{theo}\label{theo_tilted}

Let $A$ be a tilted algebra, $H$ be a hereditary algebra and $T_H$ be a tilting module such that $A = \End T$. Denote by $S\subseteq (Q_H)_0$ the set of all sink in $Q_H$, $R = \{i \in (Q_H)_0 \mid P_H(i) \in \add T\}$  and let $I = \bigoplus_{r \in R} I_H(r)$. If
\begin{enumerate}
    \item[(1)] $S \subseteq R$ ;
    \item[(2)] $\cogen \tau_H T = \add \tau_H T$;
    \item[(3)]  for any  $i \notin R$, if $\varphi \colon T_o \oplus I_0 \to I_H(i)$ is a right minimal $\add(T \oplus I)$-approximation of $I_H(i)$ then $\mathrm{Ker} \varphi \in \add T$.
\end{enumerate}
then $A$ is a representation-hereditary algebra.
\end{theo}

\proof
By (1) and by Lemma \ref{slice_L} we have that $\Hom_H(T, DH)$ is the $\Ext$-injective module in $\La_A$. And by (2) and by Lemma \ref{cogen_tau} we have $\ind A \setminus \La_A \subseteq \add DA$.
If $X_A \in \ind A$ is a non projective and non injective module such that $\Hom_A(DA, X) =0$ then $X \in \La_A$ and the projective cover $\pi \colon P_0(X) \to X$ is a  minimal right $\add(A \oplus DA)$-approximation of $X$ and, as $ \pd X = 1$, we have  $\mathrm{Ker} \pi$ projective.
If $X_A \in \ind A$ is a non projective and non injective module such that $\Hom_A(DA, X) \neq 0$ then $X_A \in \La_A$ is an Ext-injective (\cite{ACT04} 3.1), that is, $X_A = \Hom_H(T, I_H(i))$ for some injective $H$-module $I_H(i)$ and $P_H(i) \notin \add T$ (\cite[VI-Lemma 4.9]{ASS}).

Let $T_0 \in \add T$, $I_0 \in \add I$ and $\varphi\colon T_0 \oplus I_0 \to I_H(i)$ a minimal right $\add(T \oplus I)$-approximation of $I_H(i)$. As $I_H(i) \in \mathcal{T}(T) = \gen T $ then $\varphi$ is an epimorphism. By (3) in hypothesis we have $\mathrm{Ker} \varphi = T_1 \in \add T$. Applying the functor $\Hom_H(T,-)$ we obtain an epimorfism
$$f \colon \Hom_H(T, T_0) \oplus \Hom_H(T,I_0) \to \Hom_H(T, I_H(i)) = X_A $$
with $\mathrm{Ker}f = \Hom(T, \mathrm{Ker} \varphi)$ projective. This morphism $f$ is in fact a right $\add(A \oplus DA)$-approximation of $X_A$. Indeed, if $\bar I_A \in \ind A$ is injective and $\Hom_A(\bar I, X) \neq 0$ then $\bar I_A$ belongs $\La_A$ and so $\bar I_A = \Hom_H(T, I_H(j)) $ for some injective $H$-module $I_H(j)$ with $P_H(j) \in \add T$. Given $h\colon \Hom_H(T, I_H(j)) \to \Hom_H(T, I_H(i))$, consider $\psi \colon I_H(j) \to I_H(i)$ such that $\Hom_H (T, \psi) = h$. As $\varphi$ is a right $\add(T \oplus I)$-approximation of $I_H(i)$ then there exists $l\colon I_H(j) \to T_0 \oplus I_0$ such that $\varphi l = \psi$. Therefore $$h = \Hom_H(T, \varphi) = \Hom_H(T, \varphi)\circ \Hom_H(T, l) = f \circ \Hom_H(T, l). $$
Then $A$ is a hereditary representation-algebra because Theorem \ref{kernel_projective}.
\theend

Examples of a tilted algebra that are  representation-hereditary and not.
\subsection{}
\begin{ex} \begin{enumerate}
\item Let $A = \mathrm{End}_{\Lambda} (T)$ be  a tilted algebra given by the quiver 
{\scriptsize
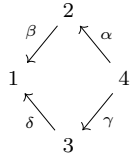
\begin{center}
\begin{tikzcd}[column sep=small]
 & 2 \ar{dl}[swap]{\beta}   &    \\
 1        &   & 4   \ar{dl}{\gamma}  \arrow{ul}[swap]{\alpha} \\
 & 3 \ar{ul}{\delta} & 
\end{tikzcd}
\end{center}
}
\noindent with the relation $\alpha \beta = 0 = \gamma \delta$. By computing the Auslander-Reiten quiver of $A$, we obtain

{\tiny
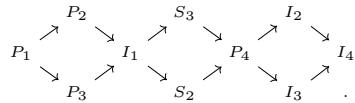
\begin{center}
\begin{tikzcd}[column sep=small, row sep=small]
& P_2 \ar{dr} &  & S_3 \ar{dr}  &  &  I_2  \ar{dr} &   \\
  P_1 \ar{ur} \ar{dr} &  & I_1 \ar{ur} \ar{dr}  & &  P_4 \ar{dr}  \ar{ur} & & I_4 \\
& P_3  \ar{ur}&  & S_2 \ar{ur}  &   & I_3 \ar{ur} & .
\end{tikzcd}
\end{center}
}

 For the simple module that are not projective and injective, we can see by the picture the left and  right minimal $A \oplus DA$-approximation with injective cokernels  and projective kernels respectively. 
  \item Let $A$ be the hereditary path algebra of the quiver $Q$

{\scriptsize
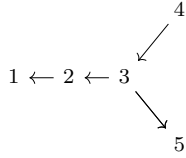
\begin{center}
\begin{tikzcd}[column sep=small]
& && 4 \ar{dl}  \\
1 & \ar{l} 2  &  3  \ar{l} \ar{dr} \ar{dr} &    \\
& && 5 
\end{tikzcd}
\end{center}
}
\noindent Consider the following tilting $A$-module $T = P(5) \oplus P(4) \oplus \tau^{-2} P(5) \oplus I(5) \oplus I(4)$, where $\tau$ is the Auslander-Reiten translation. The algebra $A = \mathrm{End}_A (T)$ is a tilted algebra given by the quiver

{\scriptsize
\begin{center}
\begin{tikzcd}
1 & \ar{l}[swap]{\delta} 2  &  3  \ar{l}[swap]{\gamma} & 4  \ar{l}[swap]{\beta}  &  5 \ar{l}[swap]{\alpha}  \\
\end{tikzcd}
\end{center}
}
\noindent subject to the relation $\alpha \beta \gamma \delta = 0$. The Auslander-Reiten quiver of $A$ is computed as

{\tiny
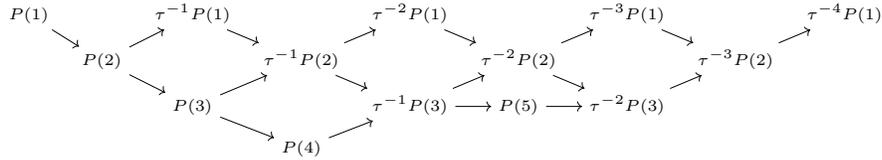
\begin{center}
\begin{tikzcd}[column sep=small, row sep=small]
P(1) \ar{dr} && \tau^{-1} P(1) \ar{dr} &  & \tau^{-2} P(1) \ar{dr}  &  &  \tau^{-3} P(1)  \ar{dr} &  & \tau^{-4} P(1) \\
 & P(2) \ar{ur} \ar{dr} &  & \tau^{-1} P(2) \ar{ur} \ar{dr}  & & \tau^{-2} P(2) \ar{dr}  \ar{ur} & & \tau^{-3} P(2) \ar{ur} & \\
&& P(3) \ar{dr}  \ar{ur}&  & \tau^{-1} P(3) \ar{ur} \ar{r} & P(5) \ar{r} & \tau^{-2} P(3) \ar{ur} & & \\
& & & P(4) \ar{ur} &  &  & & & 

\end{tikzcd}
\end{center}
}
Note that a right $A \oplus DA$-approximation of $\tau^{-4} P(1)$ does not have a projective kernel.
\end{enumerate}
\end{ex}

\vspace{1cm}
\noindent {\bf Acknowledgments.}
The first named author was partially supported by DMAT-Universidade Federal do Paran\'a and CNPq Universal Grant-405540/2023-0.
 The second named author was partially supported by Unioeste and PNPD-CAPES and visited the Universidade Federal do Paran\'a and Universidad Nacional de Mar del Plata. This study was financed in part by the Coordena\c{c}\~ao de Aperfei\c{c}oamento de Pessoal de N\'ivel Superior-Brasil(CAPES)-Finance Code 001.  The third author visited the Universidade Federal do Paran\'a.  She is supported by PICT-2021 01154 ANPCyT. 
The last author had a stay in the CEMIM-Universidad Nacional de Mar del Plata,  supported by a  posdoctoral scholarship of   CONICET - Argentina. 

\newpage



\bibliographystyle{amsplain}
\bibliography{bibliografiaHeily.bib}

\end{document}